\documentclass{amsart}

\usepackage{amsmath,amsthm,amsfonts,latexsym,amssymb,enumerate,color}
\usepackage{graphicx}
\usepackage{mathrsfs}

 \newtheorem{thm}{Theorem}[section]
 \newtheorem{cor}[thm]{Corollary}
 \newtheorem{lem}[thm]{Lemma}
 \newtheorem{prop}[thm]{Proposition}
 \theoremstyle{definition}
 \newtheorem{defn}[thm]{Definition}
 \theoremstyle{remark}
 \newtheorem{rem}[thm]{Remark}
 \newtheorem{ex}[thm]{Example}
 \numberwithin{equation}{section}

\def\CC{\mathbb C}

\def\DD{\mathbb D}
\def\GD{\mathcal D}

\def\GE{\mathcal E}
\def\GF{\mathcal F}
\def\G{\mathcal G}
\def\GH{\mathcal H}

\def\GK{\mathcal K}
\def\GL{\mathcal L}
\def\GM{\mathcal M}
\def\GN{\mathcal N}

\def\RR{\mathbb R}
\def\GS{\mathcal S}
\def\TT{\mathbb T}

\def\ker{\mathop{\rm ker}\nolimits}

\def\spam{\mathop{\rm span}\nolimits}  

\numberwithin{equation}{section}

\def\beginpf{\begin{proof}}
\def\endpf{\end{proof}}
\def\beq{\begin{equation}}
\def\eeq{\end{equation}}
\def\ds{\displaystyle}
\def\ol{\overline}
\def\til{\tilde}
\def\wt{\widetilde}

\def\sab{S_{a,b}}

\def\k{{\rm ker}}

\def\H2b{\overline{H^2_0}}

\def\PV{{\rm PV}}
\def\CP{{\rm CP}}

\def\q{\quad}
\def\th{\theta}

\begin{document}

%
%
%
%
%
%
%
%
%

\title{Paired kernels and truncated Toeplitz operators}

\author{M.~Cristina C\^amara}
\address{Center for Mathematical Analysis, Geometry and Dynamical Systems,
Instituto Superior T\'ecnico, Universidade de Lisboa, 
Av. Rovisco Pais, 1049-001 Lisboa, Portugal}
\email{cristina.camara@tecnico.ulisboa.pt} 

\author{Jonathan R.~Partington}
\address{School of Mathematics, University of Leeds, Leeds LS2~9JT, U.K.}
\email{j.r.partington@leeds.ac.uk}



\date{}

\begin{abstract} 
This paper considers paired operators in the context of the Lebesgue
Hilbert space $L^2$ on the unit circle and its subspace, the Hardy space $H^2$.
The kernels of such operators, together with
their analytic projections, which are generalizations of Toeplitz kernels, are studied. 
Inclusion relations between such kernels are considered in detail, and
the results are applied to describing the kernels of finite-rank asymmetric
truncated Toeplitz operators.
\end{abstract}

\maketitle

\centerline{\em To Ilya Spitkovsky on his 70th birthday}

\section{Introduction}\label{sec:1}

If $X$ is a Banach space, $P \in \GL(X)$ is a projection,   $Q=I-P$
is its complementary projection, and $A$ and $B$ are bounded operators on $X$, then
we use the term {\em paired operators\/} for operators of the form
\[
S_{A,B}=AP+BQ  .
\]

Paired operators appeared in the context of singular integral operators on $L^2(\RR)$ of
the form
\beq\label{eq:1jun19}
a(t)f(t) + b(t) \frac{1}{\pi i} \PV \int_\RR \frac{f(y)}{y-t} \, dy = g(t) ,
\eeq
where $a$ and $b$ are bounded functions and $g$ is given.
They were introduced by Shinbrot \cite{shinbrot}, where \eqref{eq:1jun19} was
placed in an abstract Hilbert space setting and a general method of
obtaining solutions to that equation was derived in that abstract form.

Here we shall take $X=L^2(\TT)$, where $\TT$ is the unit circle, $P$ and $Q$
are the orthogonal projections from $L^2(\TT)$ onto the Hardy space $H^2_+:=H^2(\DD)$
and its orthogonal complement $H^2_-:= \ol{H^2_0}$, respectively (we write $P=P^+$
and $Q=P^-$), and $A,B$ are multiplication operators $A=M_a$ and $B=M_b$ with $a,b \in L^\infty:=L^\infty(\TT)$. Thus we write
\beq\label{eq:2jun20}
S_{a,b}=aP^++bP^- \quad \hbox{and} \quad \Sigma_{a,b}=P^+aI+P^-bI.
\eeq
In each case the pair $(a,b)$ is called the {\em symbol pair}. These operators have been
studied by several mathematicians \cite{CGP23,CP24,DDS24,GK1,GK2,gu,MP,pross,speck}, 
often in the 
context of singular integral operators. If the operators $A$ and $B$ are invertible in $\GL(X)$, then
the paired operators are said to be of {\em normal type\/} \cite{MP}, and this is the most studied case by far.

Recently, several new results were  obtained for operators of the form 
\eqref{eq:2jun20} with $a,b \in L^\infty$ satisfying weaker conditions \cite{CGP23,CP24}:
\beq\label{eq:rev1.3}
a \ne 0 \hbox{ a.e.\ on } \TT, \quad b \ne 0 \hbox{ a.e.\ on } \TT
\eeq
and
\beq\label{eq:rev1.4}
\hbox{either} \quad a = b \quad \hbox {or} \quad a-b \ne 0 \hbox{ a.e.\ on } \TT,
\eeq
which will be standing assumptions in this paper. If $a \ne b$ we say that $(a,b)$ is {\em nondegenerate}.

Here we focus on paired operators of the form $S_{a,b}$ and
their kernels. Paired operators are
closely connected with Toeplitz operators. Indeed, on the one hand they are dilations of
Toeplitz operators; on the other hand, if $b \in \G L^\infty$ (that is, if $b$ is invertible in $L^\infty$), then $S_{a,b}$ is equivalent after extension
to the Toeplitz operator $T_{a/b}$ \cite{C17}, implying that those two types of operator share many properties such as invertibility and
Fredholmness, and their kernels are isomorphic.
(See, for instance, \cite{BTsk} for the definition and properties of operators that are equivalent after extension.)

Paired kernels, i.e., kernels of paired operators, are thus particularly connected with Toeplitz kernels, i.e., kernels of bounded Toeplitz operators. If $b \in \G L^\infty$, not only do we have
\beq
\ker S_{a,b} \cong \ker T_{a/b}
\eeq
but in fact
\beq
P^+ \ker S_{a,b} = \ker T_{a/b}.
\eeq
The spaces of the form $P^+ \ker S_{a,b}$, denoted by $\k^+ S_{a,b}$, are called
{\em projected paired kernels}. They appear in different forms: in the study of scalar-type kernels of block Toeplitz operators \cite{CP20}, in the study of nearly invariant subspaces for shift semigroups \cite{LP22}, in the characterization of the ranges of finite rank truncated Toeplitz operators (TTO)
\cite{CC19}, and so on. In Section~\ref{sec:2} we review some of the basic properties
of such projected paired kernels, and   describe their properties related to
near invariance, a property that Toeplitz kernels are known to possess.
From this point of view it
 is   natural to ask whether further known properties of Toeplitz kernels
can be extended, or have analogues, for paired kernels. 
In Section~\ref{sec:3} we focus   on the relations between Toeplitz kernels
corresponding to Toeplitz operators whose symbols are related by a multiplication operator,
such as the known relation
\beq\label{eq:1.7jun20}
\ker T_{\theta g } \subseteq \ker T_g \quad \hbox{if } \theta \hbox{ is an inner function},
\eeq
where the inclusion is strict if $\theta$ is not a constant \cite{CMP16,  CP14}; moreover,
\beq h_+ \ker T_{h_+g} \subseteq \ker T_g
\quad \hbox{if } h_+ \in H^\infty,
\eeq
where the inclusion is an equality if $h_+ \in \G H^\infty$; also
\beq\label{eq:1.9jun20}
\ker T_{h_- g} \supseteq \ker T_g \quad \hbox{if } h_- \in \ol{H^\infty},
\eeq
where the inclusion is an equality if $\ol {h_-}$ is outer.

The relations above can be roughly seen as meaning that, when we multiply the symbol of a Toeplitz operator by a given function, the kernel may become ``smaller'' or ``larger''
-- which cannot be expressed otherwise unless we are dealing with finite-dimensional kernels.

We wish to study and express the relations between paired kernels, corresponding to symbol
pairs related by multiplication operators, similarly.
It is clear that
\beq \ker S_{a,b}= \ker S_{ah,bh}  \quad \hbox{ for all } h \in L^\infty   \hbox{ with }   h \ne 0 \hbox{ a.e.}, 
\eeq
but inclusion relations analogous to \eqref{eq:1.7jun20}--\eqref{eq:1.9jun20}, which in
the finite-dimension\-al case express the fact that the dimension of one kernel is smaller
than that of the other, are impossible to obtain because the
only paired kernels contained in $\ker S_{a,b}$ are $\{0\}$ and $\ker S_{a,b}$ itself
\cite{CGP23}.

It can be shown, however, that under our standing assumptions \eqref{eq:rev1.3} and \eqref{eq:rev1.4}
 the projection $P^+$ maps $\ker S_{a,b}$ bijectively onto
$P^+ \ker S_{a,b}= \k^+ S_{a,b}$,
indeed $P^+\phi= \dfrac{b}{b-a}\phi$, given that \eqref{eq:rev1.4} holds
 (see \cite[Prop. 2.1]{CP24} for further details).

Thus we
study here the relations between the projected paired kernels with symbol pairs related by multiplication operators,
for which one can obtain inclusions and equalities analogous to \eqref{eq:1.7jun20}--\eqref{eq:1.9jun20}. These relations have also been sudied in \cite{CP24}, but we present them here in Section \ref{sec:3} for completeness, on the one hand, and also because several results are
generalized. In Section \ref{sec:4} we apply them to study and
describe the kernels of a class of asymmetric truncated Toeplitz operators (ATTO)
of the form
\beq\label{eq:11jul1.11}
 A^{\theta,\alpha}_\phi: K_\theta \to K_\alpha, \qquad A^{\theta,\alpha}_\phi=P_\alpha \phi P_\theta{}_{|K_\theta},
\eeq
where $\theta,\alpha$ are inner functions, $K_\theta,K_\alpha$ the associated model spaces,
and $P_\theta,P_\alpha$ the orthogonal projections from $L^2$ onto these model spaces \cite{CP17}.
In addition we take $\phi \in L^\infty$ with
\beq \phi=\ol\theta \frac{A_{1-}}{B_{1-}}-\alpha \frac{A_{2+}}{B_{2_+}},
\qquad A_{1-},B_{1-} \in \ol{H^\infty}, \quad A_{2+},B_{2+} \in H^\infty,
\eeq
such that $(A_{2+},B_{2+}) \in \CP^+$ and $(A_{1-},B_{1-}) \in \CP^-$
(corona pairs, defined  in Section  \ref{sec:4}),
which are motivated by the study of ATTO of finite rank. 
This is done by resorting to a characterization of certain kernels of block
Toeplitz operators that are called scalar-type kernels \cite[Thm.~3.1]{CP20}
and using the results of the previous section: this involves establishing certain
relations between the corona theorem and Hankel operators.
It is shown, in particular, that $\ker A^{\theta,\alpha}_\phi$
does not depend on $\alpha$ nor on the numerators $A_{1-}$ and $A_{2+}$
if $\dim K_\alpha$ is siufficiently large, which may be seen as rather a surprising result.

\section{Projected paired kernels and near invariance} \label{sec:2}

In this section we summarize some of the basic properties of projected paired kernels. 
A more detailed treatment is given in \cite{CP24}.
We remind the reader of our standing assumptions \eqref{eq:rev1.3} and \eqref{eq:rev1.4}.

The projected paired kernel $\k^+ \sab =P^+ \ker S_{a,b}\subset H^2_+$ was introduced in Section \ref{sec:1}, and analogously we define $\k^-\sab = P^- \ker S_{a,b} \subset H^2_-$.
Note that $\phi_- \in \k^- S_{a,b}$ if and only if for some $\phi_+ \in H^2_+$
we have $a \phi_+ + b \phi_ - = 0$; that is, $ \ol a \ol z\ol{\phi_+} + \ol b\ol z \ol{\phi_-}=0$,
and so we see that  
\beq\label{eq:9jul2.1}
\k^- S_{a,b} = \ol z \ol{\k^+ S_{\ol b,\ol a}}.
\eeq
We may rewrite this by saying that the antilinear mapping
\beq\label{eq:10jul2.1A}
T: \k^+ S_{\ol b,\ol a}\to  \k^- S_{a,b}, \qquad \psi_+ \mapsto \ol z \ol{\psi_+}
\eeq
is a bijection.

There is likewise a bijection 
\beq\label{eq:9jul2.3}
M_{-a/b}: \k^+ S_{a,b} \to \k^- S_{a,b}, \qquad \phi_+ \mapsto -\frac{a}{b} \phi_+
\eeq
(this is well-defined under the standing assumption \eqref{eq:rev1.3}).
In general we omit the index $a/b$ and write simply $M=M_{-a/b}$. So for $\phi_+ \in H^2_+$ we have 
\beq
\phi_- = M \phi_+ \iff \phi_- \in H^2_- \quad \hbox{and} \quad a\phi_++b\phi_- = 0.
\eeq

Although $\ker S_{a,b}$ is naturally a closed subspace of $L^2$, its projections
$\k^+ S_{a,b}$ and $\k^- S_{a,b}$ need not be. Further details are given in \cite{CP24}.

We now proceed to a discussion of near invariance.

\begin{defn}\label{def:15jul2.1}
Let $\GM$ be a subspace of $H^2_+$. Then we say that $\GM$ is {\em nearly invariant\/} 
or {\em nearly $S^*$-invariant\/} if and only if
\beq\label{eq:nivdef}
f \in \GM, \ f(0)=0 \implies S^*f \in \GM.
\eeq
\end{defn}
We see that \eqref{eq:nivdef} is equivalent to the condition
\[
f \in \GM, \quad \ol z f \in H^2_+ \implies \ol z f \in \GM.
\]

It is well known that kernels of Toeplitz operators   are nearly invariant, since
if $\phi \in \ker T_g$, so that
$g \phi=\psi_-$ for some $\psi_- \in H^2_-$, and $\ol z \phi \in H^2_+$ we have
$g \ol z \phi = \ol z \psi_- \in H^2_-$ and so $\ol z \phi \in \ker T_g$.

Similarly $\k^+ S_{a,b}$, is nearly invariant, since if 
$\phi_+ \in \k^+ S_{a,b}$ then there is a $\phi_- \in H^2_-$ such that
$a\phi_+ + b \phi_-=0$. Now if $\ol z\phi_+$ in $H^2_+$, then the identity
$a (\ol z \phi_+) + b \ol z \phi_- =0$ shows that $\ol z \phi_+ \in \k^+ S_{a,b}$,
since $\ol z\phi_- \in H^2_-$.  

Now, although $\k^+ S_{a,b}$ need not be a closed subspace
(see Example \ref{ex:9jul4.10} below), it is known
that the closure of a nearly invariant subspace is still nearly invariant
\cite[Prop. 4.4]{CP24}. Closed nearly invariant spaces have been considered by
 many mathematicians,
beginning with the work of Hitt~\cite{hitt} and Hayashi~\cite{hayashi}.

Let us now broaden the discussion by letting $\theta   \in H^\infty$   be a
non-constant  inner function. 

\begin{defn}
A subspace $\GM \subset H^2_+$ is {\em nearly $\ol\theta$-invariant\/}
if $f \in \GM$, $\ol \theta  f\in H^2_+$ imply that
$\ol\theta f \in \GM$.
\end{defn}

An argument similar to that above, replacing $\ol z$ by $\ol \th$
shows that projected paired kernels $\k^+ S_{a,b}$ are nearly $\ol\th$-invariant.

Note that $\phi_+ \in \k^+ S_{a,b} \cap \theta H^2_+$ if and only if we can write
$\phi_+ = \theta \psi_+$ with $\psi_+ \in H^2_+$ and $a \theta \psi_+ + b \psi_- =0$ for
some $\psi_- \in H^2_-$. Thus we conclude that
\beq
\k^+ S_{a,b} \cap \theta H^2_+ = \theta\, \k^+ S_{a\theta,b}.
\eeq
Equivalently, we may write this as
\beq
\phi_+ \in \k^+ S_{a,b}, \quad \ol\theta \phi_+ \in H^2_+ \iff \phi_+ \in \theta\, \k^+ S_{a\theta,b},
\eeq
and, thus, saying that $\k^+ S_{a,b}$ is nearly $\ol\theta$-invariant is equivalent to saying that
\[
\theta\, \k^+ S_{a\theta, b} \subset \k^+ S_{a,b}.
\]
By similar arguments we have also
\beq\label{eq:9jul3.3}
\k^- S_{a,b} \cap \ol\theta H^2_- = \ol\theta\, \k^- S_{a,\ol\theta b},
\eeq
and so
\beq
\phi_- \in \k^- S_{a,b}, \quad \theta \phi_- \in H^2_- \iff \phi_- \in \ol\theta \,\k^- S_{a,\ol\theta b}.
\eeq
In particular, 
\beq \k^+ S_{a,b} \cap zH^2_+ = z \k^+ S_{az,b},
\eeq
consisting of the elements of $\k^+ S_{a,b}$ that vanish at $0$; analogously, $\ol z \k^- S_{a,\ol z b}$ consists of the functions $\phi_- \in \k^- S_{a,b}$ such that $z \phi_-$ vanishes at $\infty$.

Also note that, for $\phi_- = M \phi_+$ (i.e., such that $a\phi_++b \phi_-=0$ as in 
\eqref{eq:9jul2.3}), we have
\beq\label{eq:9jul3.6}
\phi_+ \in \theta\, \k^+ S_{a \theta,b} \iff \phi_- \in \k^- S_{a,\ol\theta b}
\eeq
and
\beq\label{eq:9jul3.7}
\phi_+ \in \k^+ S_{a\theta,b} \iff \phi_- \in \ol\theta\, \k^- S_{a,\ol\theta b}.
\eeq

\section{Inclusion relations}
\label{sec:3}

Let us now study the relations between two paired kernels with
symbol pairs that are related by multiplication operators, by means of certain inclusion
relations between their projections into $H^2_\pm$. Although some of these relations were obtained in \cite{CP24} we present here all proofs for completeness, which also allows us in some cases to present alternative proofs. We again remind the reader of our standing assumptions \eqref{eq:rev1.3} and \eqref{eq:rev1.4}.

\begin{prop}\label{prop:9jul4.1}
Let $\theta_1$ and $\theta_2$ be inner functions. Then 
\beq \label{eq:9jul4.1}
\k^+ S_{a\theta_1,b \ol\theta_2} \subseteq \k^+ S_{a,b},
\eeq
and, for $\k^+ S_{a,b} \ne \{0\}$, the inclusion is proper if and only if 
$\theta_1 \theta_2 \not\in \CC$; that is,  
$\theta_1$ or
$ \theta_2$ is non-constant.
\end{prop}
\beginpf
Let $\theta=\theta_1\theta_2$; then $\k^+ S_{a\theta_1,b\ol{\theta_2}}=\k^+ S_{a\theta,b}$. Let $\phi_+ \in \k^+ S_{a\theta,b}$, so that
$a\theta\phi_+ + b \psi_- =0$ for some $\psi_- \in H^2_-$. Then we have
$a\phi_+ + b (\ol\theta \psi_-)=0$ and so $\phi_+ \in \k^+ S_{a,b}$.

If $\theta \in \CC$, then the inclusion \eqref{eq:9jul4.1} is an equality.
Otherwise, if $\theta\not\in \CC$, then there exists $\phi_- \in \k^- S_{a,b}$ such that $\phi_- \not\in \ol\theta H^2_-$, by near invariance: see Propositions 4.1 and 4.5 in \cite{CP24}.

Since $\phi_- \in \k^- S_{a,b}$ but $\phi_- \not\in \ol\theta \k^- S_{a,\ol\theta b}$ by
\eqref{eq:9jul3.3}, for $\phi_+ = M^{-1} \phi_-$ we must have that
$\phi_+ \in \k^+ S_{a,b}$ but $\phi_+ \not\in \k^+ S_{a\theta,b}$, by \eqref{eq:9jul3.7},
so the inclusion is strict.
\endpf

By \eqref{eq:9jul2.1} we have the following.

\begin{cor}\label{cor:10jul4.2}
Let $\theta_1$ and $\theta_2$ be inner functions. Then
\beq
\k^- S_{a\theta_1,b\ol\theta_2} \subseteq \k^- S_{a,b},
\eeq
and, for $\k^- S_{a,b} \ne \{0\}$, the inclusion is strict if and only if $\theta_1\theta_2 \not\in \CC$.
\end{cor}

\begin{prop}\label{prop:9jul4.3}
Let $h_- \in \ol{H^\infty}$. Then we have:\\
(i) $\k^+ S_{a,bh_-} \subseteq \k^+ S_{a,b} \subseteq \k^+ S_{ah_-,b}$;\\
(ii) If $\ol{h_-}$ is outer then
\beq
\k^+ S_{a,bh_-} = \k^+ S_{a,b} \quad \hbox{if} \quad  \frac{a}{bh_-} \in L^\infty,
\eeq
and
\beq\label{eq:9jul4.4}
\k^+ S_{ah_-,b}=\k^+ S_{a,b} = \ker T_{a/b} \quad \hbox{if} \quad a/b \in L^\infty;
\eeq
(iii) If $\ol h_-$ is not outer, then the inclusions in (i) are strict.
\end{prop}
\beginpf
The inclusions in (i) are immediate consequences of the definition of a projected paired
kernel. Regarding (ii), we have that if
$a \phi_+ + b\phi_- =0$ with $\phi_\pm \in H^2_\pm$, then
\[
a\phi_+ + bh_-(h_-^{-1}\phi_-)=0 \quad \hbox{with} \quad h_-^{-1} \phi_-=-\frac{a}{bh_-}\phi_+ \in L^2 \cap \ol z \ol{\GN_+},
\]
where $\GN_+$ denotes the Smirnov class (see, e.g.,  \cite{Nik}), if $\dfrac{a}{bh_-} \in L^\infty$. 
So $h_-^{-1}\phi_- \in H^2_-$, and therefore $\phi_+ \in \k^+ S_{a,bh_-}$.

This, together with the first inclusion in (i), proves the first equality in (ii), and the second equality is proved analogously.

As for (iii), if $\ol{h_-}$ has a non-constant inner factor $\theta$, then it follows from Prop.~\ref{prop:9jul4.1} and from (i) that, if $h_- = \ol\theta \tilde h_-$ with $\tilde h_- \in \ol{H^\infty}$, then
\[
\k^- S_{a,b\ol\theta \tilde h_-} \subsetneq \k^+ S_{a,b \tilde h_- }\subseteq \k^+ S_{a,b},
\]
and analogously for the second inclusion.
\end{proof}

\begin{rem}\label{rem:9jul4.4}{\rm
Note that the inclusions in Prop.~\ref{prop:9jul4.3}~(i)
may be strict even if $\ol {h_-}$ is outer: see Example~\ref{ex:9jul4.10} below.
Also note that from \eqref{eq:9jul4.4} we have that, if $g \in L^\infty$, $\ol{h_-}$ is outer
and $h_+ \in \G H^\infty$, then
\beq\label{eq:9jul4.5}
\k^+ S_{gh_-,h_+} = \ker T_{gh_+^{-1}} = h_+ \ker T_g.
\eeq
}
\end{rem}

\begin{cor} Let $z_0 \in \TT$. Then\\
(i) if $a/b \in L^\infty$, then $\k^+ S_{a(z-z_0),b} = \k^+ S_{az,b}$;\\
(ii) if  $\dfrac{a}{b(z-z_0)} \in L^\infty$, then $\k^+ S_{a,b(z-z_0)} = \k^+ S_{a,bz}$.
\end{cor}

\beginpf
Take $h_- = (z-z_0)/z \in \ol{H^\infty}$, which is such that $\ol{h_-}$ 
is outer in $H^\infty$. Then (i) and (ii) follow from the second equality and first equality in
Prop.~\ref{prop:9jul4.3}~(ii), respectively.
\endpf
Note that
\beq\label{eq:12jul4.6}
\k^+ S_{a,b} = \ker T_{a/b}, \quad \hbox{if} \quad a/b \in L^\infty,
\eeq
so we have
\[
\ker T_{g(z-z_0)}=\ker T_{g z} \quad \hbox{for} \quad z_0 \in \TT \quad \hbox{and} \quad g \in L^\infty.
\]
Taking \eqref{eq:9jul2.1} into account, we also have:
\begin{cor} Let $h_+ \in H^\infty$; then
\[
\k^- S_{ah_+,b} \subseteq \k^- S_{a,b} \subseteq \k^- S_{a,bh_+}.
\]
\end{cor}

\begin{prop}\label{prop:9jul4.7}
Let $h_+ \in H^\infty$. Then we have:\\
(i) $h_+ \k^+ S_{ah_+,b} \subseteq \k^+ S_{a,b}$ and
$h_+ \k^+ S_{a,b} \subseteq \k^+ S_{a,bh_+}$;\\
(ii) if $h_+ \in \G H^\infty$, or if $h_+$ is outer with $\dfrac{b}{ah_+} \in L^\infty$, then
\[
h_+ \k^+ S_{ah_+,b} = \k^+ S_{a,b}   ;
\]
(iii)
if $h_+ \in \G H^\infty$, or if $h_+$ is outer with $\dfrac{b}{a} \in L^\infty$, then
\[
h_+ \k^+ S_{a,b} = \k^+ S_{a,bh_+};
\]
(iv) if $h_+$ is not outer, then the inclusions in (i) are strict.
\end{prop}

\beginpf
(i) The first inclusion is obvious, since for   $\phi_- \in H^2_-$
\[
(ah_+) \phi_+ + b \phi_-=0  \iff a(h_+ \phi_+) + b \phi_- = 0.
\]
Regarding the second inclusion, if $\phi_+ \in \k^+ S_{a,b}$ then
$a(h_+ \phi_+) + (b h_+) \phi_- =0$, so $h_+ \phi_+ \in \k^+ S_{a,bh_+}$.\\
(ii) It is left to show that $\k^+ S_{a,b} \subseteq h_+ \k^+ S_{ah_+,b}$. We have that
\[
a\phi_+ + b\phi_-=0 \iff (ah_+) \frac{\phi_+}{h_+} + b \phi_- = 0 \iff \frac{\phi_+}{h_+}=-\frac{b}{ah_+} \phi_-.
\]
If $h_+ \in \G H^\infty$, then ${\phi_+}/{h_+} \in H^2_+$; if $h_+$ is outer and
$\dfrac{b}{ah_+} \in L^\infty$, then 
$\phi_+/h_+ \in L^2 \cap \GN_+$ and so $\phi_+/h_+ \in H^2_+$.
It follows that $\phi_+/h_+ \in \k^+ S_{ah_+,b}$.\\
(iii) Now it remains to show that
$\k^+ S_{a,bh_+} \subseteq h_+ \k^+ S_{a,b}$.
We have
\[ a\phi_+ + (bh_+)\phi_- = 0 \iff a \phi_+/h_+ + b \phi_-=0. \] 
If $h_+ \in \G H^\infty$, then $\phi_+/h_+ \in H^2_+$; if $h_+$ is outer and
$b/a \in L^\infty$ then
$\phi_+/h_+ \in L^2 \cap \GN_+ = H^2_+$, so $\phi_+ \in h_+ \k^+ S_{a,b}$.
\\
(iv) If $h_+$ is not outer, we cannot have $\k^+ S_{a,b} \subseteq h_+ \k^+ S_{ah_+,b} \subseteq h_+ H^2_+$, nor $\k^+ S_{a,bh_+} \subseteq h_+ \k^+ S_{a,b} \subseteq h_+ H^2_+$, by
near invariance (Cor.~4.2 in \cite{CP24}).
\endpf

\begin{rem}{\rm
Note that the inclusion in (i) above may be strict even if $h_+$ is outer (see also Remark \ref{rem:9jul4.4}).

Recalling 
\eqref{eq:12jul4.6},
it follows from the second relation in Prop.~\ref{prop:9jul4.7}~(i) that
\beq\label{eq:9jul4.6}
h_+ \ker T_a \subseteq \k^+ S_{a,h_+}
\eeq
for any $h_+ \in H^\infty$. In particular, since $\ker T_a \ne \{0\}$ when $a \in \ol{H^\infty}$ 
and $\ol a$ is not outer, it follows from \eqref{eq:9jul4.6}
that if $h_+$ or $\ol{h_-}$ is not outer, then
\beq\label{eq:23aug3.8}
\k^+ S_{h_-,h_+} \ne \{0\}.
\eeq
}
\end{rem}

\begin{cor}
Let $h_- \in \ol{H^\infty}$. Then
\beq
h_- \k^- S_{a,bh_-} \subseteq \k^- S_{a,b} \quad \hbox{and} \quad
h_- \k^- S_{a,b} \subseteq \k^- S_{ah_-,b}.
\eeq
These inclusions are equalities if $\ol{h_-} \in \G H^\infty$ or
$\ol{h_-}$ is outer, $\dfrac{a}{bh_-} \in L^\infty$; they are strict if $\ol{h_-}$
is not outer.
\end{cor}

In the following example we illustrate some of the inclusion relations presented above.

\begin{ex}\label{ex:9jul4.10}{\rm
Let $\theta$ be the singular inner function defined by
\beq 
\theta(z) = \exp \left( \frac{z-1}{z+1} \right), \qquad (z \in \TT).
\eeq
Then, by Prop.~\ref{prop:9jul4.7}~(ii),
\beq
(z+1)K_\theta = (z+1) \k^+ S_{\ol\theta,1}=\k^+ S_{\ol\theta,z+1}.
\eeq
By Prop.~\ref{prop:9jul4.3}~(i) and \eqref{eq:9jul4.5},
\beq
(z+1)K_\theta = \k^+ S_{\ol\theta \ol z,(z+1)/z} \subseteq \k^+ S_{\ol\theta \ol z,1}=\ker T_{\ol\theta \ol z}=K_{\theta z},
\eeq
where the inclusion is strict because, as shown in \cite{LP22}, $(z+1)K_\theta$ is not a closed subspace so, in particular, it cannot be a Toeplitz kernel. We thus see that the inclusion in 
Prop.~\ref{prop:9jul4.3}~(i) may be strict if $\ol {h_-}$ is outer in $H^\infty$ but the condition
$\dfrac{a}{b h_-} \in L^\infty$ is not satisfied, as happens in this case.
}
\end{ex}

\begin{rem}\label{rem:9jul4.11}{\rm
Note that any space $uK_\theta \subset H^2_+$ with $u \in H^2_+$ outer, as in the example above, can be expressed as $\k^+ S_{a,b}$ if we allow for symbol pairs in $(L^2)^2$, with the domain
$\GD= \{ \phi \in L^2: a \phi_+ + b \phi_- \in L^2 \}$.
Indeed, following the same reasoning as in Prop.~\ref{prop:9jul4.7}~(ii), we have
that $uK_\theta = \k^+ S_{\ol\theta,u}$.

On the other hand, if $u$ is not outer, then $uK_\theta$ cannot be expressed as a projected paired kernel, even if we allow
for $L^2$ symbols, since this space cannot be contained in $\alpha H^2_+$ if $\alpha$ is a non-constant inner function, by
an argument similar to that of Corollary~4.2 in \cite{CP24}.
}
\end{rem}

We shall be interested in studying projected paired kernels of the form 
\beq\label{eq:9jul4.10a}
\k^+ S_{\alpha B_-,B_+} \quad \hbox{with} \quad B_+,\ol{B_-} \in H^\infty.
\eeq
These are strictly contained in $\k^+ S_{B_-,B_+}$, where we assume that either $\ol{B_-}$ or $B_+$ is not outer, so that 
$\k^+ S_{B_-,B_+} \ne \{0\}$ by \eqref{eq:23aug3.8}. In particular, we shall be interested in establishing conditions for
\eqref{eq:9jul4.10a} to be $\{0\}$.

It is thus convenient to see what form the previous inclusions take when $h_+$ or $\ol{h_-}$ are inner functions.

\begin{prop}\label{prop:9jul4.12}
Let $\theta_1$ and $\theta$ be non-constant inner functions such that 
$\theta_1 \prec \theta$ (that is,
$\theta_1$ properly divides $\theta$). Then,
\\
(i) $\k^+ S_{a\theta,b} \subsetneq \k^+ S_{a\theta_1,b} \subsetneq \k^+ S_{a,b}
\subsetneq \k^+ S_{a\ol{\theta_1},b} \subsetneq \k^+ S_{a \ol\theta,b}$;
\\
(ii) $\theta \k^+ S_{a\theta,b} \subsetneq \theta_1 \k^+ S_{a\theta_1,b} \subsetneq \k^+ S_{a,b}$,\\
assuming that the right-hand side of each inclusion is different from $\{0\}$.
\end{prop} 
\beginpf
(i) follows from Prop.~\ref{prop:9jul4.1} and (ii) follows from Prop.~\ref{prop:9jul4.7}.
\endpf

Analogously, we have:

\begin{prop}
Let $\theta_1$ and $\theta$ be non-constant inner functions such that $\theta_1 \prec \theta $. Then\\
(i) $\k^- S_{a,\ol\theta b} \subsetneq \k^- S_{a, \ol{\theta_1}b}
\subsetneq \k^- S_{a,b} \subsetneq
\k^- S_{a,b\theta_1} \subset \k^- S_{a,\theta b}$;\\
(ii) $\ol\theta \k^- S_{a,\ol\theta b} \subsetneq \ol{\theta_1} \k^- S_{a,{\ol\theta_1}b} \subsetneq \k^- S_{a,b}$,\\
assuming that the right-hand side of each inclusion is different from $\{0\}$.
\end{prop}

These inclusion relations allow us to answer practically, in some sense, the question of how much ``smaller'' $\k^+ S_{a\theta,b}$ or $\theta \k^+ S_{a\theta,b}$
are than $\k^+ S_{a,b}$ -- a question which naturally arises from the strict inclusions of Proposition~\ref{prop:9jul4.1}, 
and, in an analogous form, from Corollary~\ref{cor:10jul4.2}, assuming that $\theta$ is a non-constant. Indeed, if
$\theta_0=\theta$ and $\theta_i$ $(i=1,2,\ldots)$ are non-constant inner functions such that
$\theta_0 \succ \theta_1 \succ \theta_2 \succ \ldots$ and $\k^+ S_{a\theta,b} \ne \{0\}$, then we have
\beq \label{eq:10jul4.12}
\k^+ S_{a\theta,b} \subsetneq \k^+ S_{a\theta_1,b} \subsetneq \k^+ S_{a\theta_2,b}
\subsetneq \ldots \subsetneq \k^+ S_{a,b},
\eeq
and
\beq \label{eq:10jul4.13}
\theta \k^+ S_{a\theta,b} \subsetneq \theta_1 \k^+ S_{a\theta_1,b} \subsetneq \theta_2 \k^+ S_{a\theta_2,b} \subsetneq \ldots \subsetneq \k^+ S_{a,b},
\eeq
which, according to \eqref{eq:9jul3.6} and \eqref{eq:9jul3.7}, are equivalent to
\beq
\ol\theta \k^- S_{a,\ol\theta b} \subsetneq \theta_1 \k^- S_{a,\ol{\theta_1}b} 
\subsetneq \theta_2 \k^- S_{a,\ol{\theta_2}b} 
\subsetneq \ldots \subseteq \k^- S_{a,b},
\eeq
and
\beq
\k^- S_{a,\ol\theta b} \subsetneq \k^- S_{a,\ol{\theta_1}b} \subsetneq
\k^- S_{a,\ol{\theta_2}b} \subsetneq \ldots \subsetneq \k^- S_{a,b}.
\eeq
If $\theta$ is a finite Blaschke product, then we have the following.

\begin{prop} \label{prop:12jul4.14}
\cite[Thm. 6.8]{CP24} 
Let $B$ be a finite Blaschke product. Then\\
(i) $\k^+ S_{aB,b}$ is finite-dimensional if and only if $\k^+ S_{a,b}$ is finite-dimensional,
and analogously for $\k^+ S_{a,Bb}$;\\
(ii) if $\dim \k^+ S_{a,b} =d < \infty$ and $\dim K_B \ge d$, then
$\k^+ S_{aB, b}=\{0\}$.\\
(iii) if $\dim \k^+ S_{a,b} = d < \infty$ and $\dim K_B = k < d$, then
$\dim \k^+ S_{aB,b}=d-k$.
\end{prop}
From   \eqref{eq:10jul4.12} and \eqref{eq:10jul4.13} we also have:

\begin{prop}
Let $B=B_0$ be a finite Blaschke product of degree $k$, and let $B_1 ,B_2, \ldots, B_k$ be finite
Blaschke products such that
$B \succ B_1 \succ B_2 \succ \ldots \succ B_k \in \CC$. Moreover, suppose that
$\dim \k^+ S_{a,b} > k$.
Then there exist $\phi_{i+} \in \k^+ S_{aB_i,b} \setminus \k^+ S_{aB_{i-1}b}$, for $i=1,2,\ldots,k$,
such that
\beq \label{eq:10jul4.16}
\k^+ S_{a,b} = \k^+ S_{aB,b} + \spam \{\phi_{1+},\ldots,\phi_{k+}\}.
\eeq
and there exist $\psi_{1+} \in B_i \k^+ S_{aB_i} \setminus B_{i-1} \k^+ S_{a,B_{i-1},b}$,
for $i=1,2,\ldots,k$, such that
\beq \label{eq:10jul4.17}
\k^+ S_{a,b} = B \k^+ S_{aB,b} + \spam\{ \psi_{1+},\ldots,\psi_{k+} \}.
\eeq
\end{prop}

\beginpf
It is enough to prove the result when $B$ is a finite Blaschke product of degree $1$.
Let us take $B=z$ for simplicity: the general argument is similar.
Since $z \k^+ S_{az,b} \subsetneq \k^+ S_{a,b}$, there exists
$\wt\psi_+ \in \k^+ S_{a,b} \setminus z \k^+ S_{az,b}$, so $\wt \psi_+(0) \ne 0$.
Now for $\phi_+ \in \k^+ S_{a,b}$ define
\[
\wt \phi_+= \wt \psi_+(0)\phi_+ - \phi_+(0) \wt \psi_+.
\]

We have that 
\[
\wt \phi_+ \in \k^+ S_{a,b} \cap zH^2_+ = z \k^+ S_{az,b}
\]
and
\[
\phi_+ = \underbrace{\frac{\wt \phi_+}{\wt \psi_+(0)}}_{\in\, z \k^+ S_{az,b}} + \underbrace{\frac{\phi_+(0)}{\wt \psi_+(0)} \wt \psi_+}_{\in\, \k^+ S_{a,b} \setminus z \k^+ S_{az,b}},
\]
so we conclude that \eqref{eq:10jul4.16}  holds for $B=z$.

By  \eqref{eq:9jul3.6} it follows from this that
\[
\k^- S_{a,b} = \k^- S_{a,\ol z b} + \spam \{ \wt \phi_- \},
\]
where $\wt \phi_-  = M \wt \phi_+ \in \k^- S_{a,b} \setminus \k^- S_{a,\ol z b}$. Now, 
using \eqref{eq:10jul2.1A} we conclude that
\[
\k^+ S_{a,b} = \k^+ S_{az,b} + \spam \{ \phi_+\},
\]
where $\phi_+ \in \k^+ S_{a,b} \setminus \k^+ S_{az,b}$.
\endpf

Finally, note that for a projected paired kernel of the form \eqref{eq:9jul4.10a}, if $B_-= \ol{\theta_1}\ol{O_1}$ with $\theta_1$ inner and $O_1$ outer in $H^\infty$, and
$B_+= \theta_2 O_2$ with $\theta_2$ inner and $O_2 \in \G H^\infty$, then by 
\eqref{eq:9jul4.5} we have that
\beq \label{eq:12jul4.18}
\k^+ S_{\alpha B_-,B_+}= \k^+ S_{\alpha \ol{\theta_1}\ol{\theta_2} \ol{O_1},O_2}= O_2 \ker T_{\alpha \ol\theta}
\eeq
with $\theta=\theta_1\theta_2$, so the question whether $\k^+ S_{\alpha B_-,B_+}=\{0\}$
reduces to the question whether $\ker T_{\alpha \ol\theta}=\{0\}$. 
This, in turn, is related to the question of existence of a multiplier from $K_\theta$ into $K_\alpha$,
studied in \cite{FR18,CP18}.
The answer to our question, in the case when one of the inner functions is
a finite Blaschke product, was given in \cite[Sec.~6]{CMP16},
but the situation is more complicated when neither inner function is a finite
Blaschke product.\\

So suppose that $\theta$ and $\alpha$ are non-constant
inner functions in $H^\infty(\DD)$ with no nontrivial common factors. The
question is: when is $\ker T_{\alpha\ol\theta}$ different from $\{0\}$?
By Coburn's lemma at least one of $\ker T_{\alpha\ol\theta }$ and $\ker T_{\theta\ol\alpha}$
will be zero.\\

Recall that the spectrum of an inner function $\sigma(\theta)$ consists of
all $\lambda \in \ol\DD$ such that $\liminf_{z \to \lambda} |\theta(z)|=0$,
and it is the union of the closure of its zero set in $\DD$ and the support of its singular measure
(if any).
Useful facts from \cite[Sec.~7.3]{GMR} are:
\begin{itemize}
\item if $\sigma(\theta) \cap \TT \ne \TT$, then $\theta$ has an
analytic continuation across $\TT \setminus \sigma(\theta)$; 
\item $\theta$ does not have an analytic continuation across any
point of $\sigma(\theta) \cap \TT$;
\item every function in $K_\theta$ has an analytic continuation across
$\TT \setminus \sigma(\theta)$.
\end{itemize}

\begin{thm}\label{thm:12jul4.16}
If $\ker T_{\alpha\ol\theta} \ne \{0\}$ we have
$\sigma(\alpha)\cap\TT \subset \sigma(\theta)\cap \TT$.
\end{thm}

\beginpf
We may suppose that 
$\sigma(\theta) \cap \TT \ne \TT$, that
$f_+ \ne 0$ and $f_+ \in \ker T_{\alpha\ol\theta }$.
That is,
$\alpha\ol\theta  f_+=f_-$ for some $f_- \in H^2_-$.
Then $f_+$ and $\alpha f_+$ both lie in $K_\theta$, and have 
analytic continuations across $\TT \setminus \sigma(\theta)$.
Hence, by division, $\alpha$ has a meromorphic continuation, but since it is inner,
this must be an analytic continuation (no poles), since it has to be unimodular a.e.\ on $\TT$.
\endpf

Note that if $\theta=\alpha\beta$ with $\beta$ inner, non-constant, then
$\ker T_{\ol\theta \alpha}=K_\beta \ne \{0\}$, and
we do indeed have 
$\sigma(\alpha)  \subset \sigma(\theta) $. Moreover, the example given in \cite[Ex.~6.3]{CMP16} corresponds to $\theta$ and $\alpha$ being
singular inner functions with  point masses at different points, so we recover
the result that $\ker T_{\ol\theta \alpha}=\{0\}$ in this case.
This leaves open many other examples, e.g., what happens if $\theta$ and $\alpha$
are Blaschke products with disjoint zero sets that both converge to $1$?

\section{Kernels of ATTO}
\label{sec:4}

We defined asymmetric truncated Toeplitz operators (ATTO) in \eqref{eq:11jul1.11} and now consider the question of describing the kernel of such an operator when the symbol has the form
\beq\label{eq:11jul5.1}
\phi = \ol\th \frac{A_{1-}}{B_{1-}} - \alpha \frac{A_{2+}}{B_{2+}} \in L^\infty,
\eeq
with $(A_{2+},B_{2+}) \in \CP^+$ and $(A_{1-},B_{1-}) \in \CP^-$. Here
$\CP^+$ denotes the set of all {\em corona pairs\/} $(f_{1+},f_{2+})$ with 
$f_{1+},f_{2+} \in H^\infty$ such that there exist $\til  f_{1+},\til f_{2+} \in H^\infty$
satisfying
\beq\label{eq:23aug4.2}
f _{1+}\til f_{1+} + f_{2+}\tilde f_{2+}=1.
\eeq
By the well-known corona theorem \cite{carleson}, there exist $\til f_{1+}$ and $\til f_{2+}$ such that
\eqref{eq:23aug4.2} holds if and only if
\beq
\inf_{z \in \DD} |f_{1+}(z)| + |f_{2+}(z)| > 0.
\eeq

$\CP^-$ is defined analogously, with $H^\infty$ replaced by $\ol{H^\infty}$.

This class of ATTO  is motivated by the study of finite-rank ATTO, which will be considered as
a particular case where $A_{1-},B_{1-}$ and $A_{2+},B_{2+}$ are polynomials (see
\cite{CP17,CP24}. The starting point of this study is a result from \cite{CP20}, which is formulated here as follows, denoting by $\GF$ the set of all complex-valued measurable functions defined almost everywhere on $\TT$.

\begin{thm}
\label{thm:11jul5.1}
 \cite[Cor. 3.4]{CP20}
Let $f=(f_1,f_2)$ and $g=(g_1,g_2)$ belong to $\GF^2$, with left inverses $\til f^T$ and $\til g^T$, respectively, where
$\til f=(\til f_1,\til f_2)$ and $\til g=(\til g_1,\til g_2)$ belong to $\GF^2$; i.e.,
\beq
\til f^T f=\til g^T g = 1.
\eeq
Suppose also that, for a given $G \in (L^\infty)^{2 \times 2}$ we have that
\beq\label{eq:11jul5.4}
Gf=g.
\eeq
Let moreover 
\beq
\GS = (\det G. f^T(H^2_+)^2) \cap (g^T(H^2_-)^2),
\eeq
and
\beq
\GK = (\til f^T \ker H_{f \til f^T}) \cap (\til g^T \ker \wt H_{g \til g^T}),
\eeq
where
\beq
\ker H_{f\til f^T} = \{ \phi_+ \in (H^2_+)^2: P^- f \til f^T \phi_+=0 \}
\eeq
and
\beq
\ker \wt H_{g \til g^T} = \{ \phi_- \in (H^2_-)^2: P^+ g \til g^T \phi_-=0 \}.
\eeq
Then, if $\GS=\{0\}$ we have that
\beq\label{eq:11jul5.9}
\ker T_G = \GK f.
\eeq 
\end{thm}
\begin{rem}{\rm
Given that $f$ is a fixed vector function and $\GK$ is a space of scalar functions we
say that $\ker T_G$ is of {\em scalar type\/} if it satisfies an equation of the form
\eqref{eq:11jul5.9}.

There are many examples of scalar-type kernels for block Toeplitz operators. We have the following
analogue of Coburn's lemma for Toeplitz operators with $2 \times 2$ symbols.
}
\end{rem}

\begin{thm}\cite{CP20}
If $G \in (L^\infty)^{2\times 2}$ and $\det G \in L^\infty \setminus \{0\}$ then either
$\ker T_G=\{0\}$ or $\ker T_G^*=\{0\}$, or
both kernels are of scalar type.
\end{thm}

Applying Theorem \ref{thm:11jul5.1} to the study of kernels of ATTO will be based on the following result.

\begin{thm}\label{thm:11jul5.4}
\cite{CP17}. Let $\th,\alpha$ be inner functions, let $\phi \in L^\infty$, and set
\beq\label{eq:11jul5.10}
G= \begin{bmatrix}
\ol\theta & 0 \\ \phi & \alpha
\end{bmatrix}.
\eeq
Denoting by $P_1$ the projection $(x,y)   \overset{P_1} \longmapsto x$,
we have that
\beq
\ker A^{\theta,\alpha}_\phi = P_1 \ker T_G.
\eeq
\end{thm}
To study the kernel of an ATTO of the form $A_\phi^{\th,\alpha}$, with $\phi$ given by \eqref{eq:11jul5.1}, we shall therefore apply Theorem~\ref{thm:11jul5.1}
to matrix symbols of the form \eqref{eq:11jul5.10}.
A particular type of triangular symbol in this class was 
also studied, in the context of almost-periodic factorization on the real line,
in the book \cite{BKS}.

It is easy to see that   \eqref{eq:11jul5.4} is satisfied for
\beq\label{eq:11jul5.12}
f= \theta \begin{bmatrix} 1 \\ \frac{A_{2+}}{B_{2+}}\end{bmatrix}, \qquad 
g= \begin{bmatrix} 1 \\ \frac{A_{1-}}{B_{1-}}\end{bmatrix},
\eeq
which have left inverses $\til f^T$ and $\til g^T$ respectively, with
\beq\label{eq:11jul5.13}
\tilde f= \begin{bmatrix}\, \ol\theta\, \\ 0 \end{bmatrix}, \qquad \til g= \begin{bmatrix}1 \\ 0 \end{bmatrix}.
\eeq
It follows from Theorems \ref{thm:11jul5.1} and \ref{thm:11jul5.4} that, taking
\beq
G = \begin{bmatrix}
\ol \theta & 0 \\ \ol\th \frac{A_{1-}}{B_{1-}} - \alpha \frac{A_{2+}}{B_{2+}} & \alpha
\end{bmatrix}
\eeq
and $f,g,\til f,\til g$ as in \eqref{eq:11jul5.12} and \eqref{eq:11jul5.13}, we have the following.

\begin{thm}\label{thm:12jul5.5}
 If
\beq\label{eq:11jul5.15}
\GS:= \left( \alpha \begin{bmatrix}1 & \frac{A_{2+}}{B_{2+}}\end{bmatrix}
(H^2_+)^2 \right)
\cap
\left ( \begin{bmatrix}1 & \frac{A_{1-}}{B_{1-}}\end{bmatrix} (H^2_-)^2 \right) = \{0\},
\eeq
then, for $\phi$ given by \eqref{eq:11jul5.1},
\beq\label{eq:12jul5.16}
\ker A^{\th,\alpha}_{\phi}=\GK \theta,
\eeq
where
\beq
\GK=  \biggl( \begin{bmatrix}\,\ol\theta & 0 \end{bmatrix} \underbrace{\ker H_{\begin{bmatrix}
1 & 0 \\ \frac{A_{2+}}{B_{2+}} & 0 \end{bmatrix}}}_{\ker H_{f\til f^T}}  \biggr)
\cap
\biggl(
\begin{bmatrix}1 & 0 \end{bmatrix}\underbrace{ \ker \wt H_{\begin{bmatrix}1 & 0 \\
\frac{A_{1-}}{B_{1-}} & 0 \end{bmatrix}}}_{\ker \wt H_{g\til g^T}} \biggr).
\eeq
\end{thm}
Let us first show that, if $K_\alpha$ is ``big'' enough, then $\GS=\{0\}$
independently from what the inner function $\th$ is.

\begin{prop}\label{prop:11jul5.6}
Let $\GS$ be defined as in \eqref{eq:11jul5.15}. Then
\beq
B_{2+} \GS \subseteq \alpha \k^+ S_{\alpha B_{1-},B_{2+}}.
\eeq
\end{prop}
\beginpf
$\GS$ consists of functions of the form
\beq
\alpha\left( \phi_{1+}+ \frac{A_{2+}}{B_{2+}} \phi_{2+}\right) \qquad \hbox{with} \q \phi_{1+},\phi_{2+} \in H^2_+,
\eeq
\beq\alpha \left( \phi_{1+}+ \frac{A_{2+}}{B_{2+}} \phi_{2+}\right) = \phi_{1-}+ \frac{A_{1-}}{B_{1-}} \phi_{2-}
\qquad \hbox{for some} \q \phi_{1-},\phi_{2-} \in H^2_-,
\eeq
so 
\begin{flalign*}
B_{1-} B_{2+} \left( \alpha(\phi_{1+}+ \frac{A_{2+}}{B_{2+}} \phi_{2+}) \right)
& = 
\alpha B_{1-} (B_{2+}\phi_{1+}+A_{2+}\phi_{2+}) \\
&= B_{2+}(B_{1-}\phi_{1-}+ A_{1-}\phi_{2-}).
\end{flalign*}
Thus 
\[
B_{1-}B_{2+} \GS \subseteq (\alpha B_{1-} H^2_+) \cap (B_{2+}H^2_- ),
\]
which, using the lemma below, means that
\[
B_{1-}B_{2+}\GS \subseteq \alpha B_{1-} \k^+ S_{\alpha B_{1-},B_{2+}} 
\]
as required.
\endpf

At this point it will be helpful  to observe that $\k^+ S_{A,-B}=\k^+ S_{A,B}$ for
all $A,B \in L^\infty$.

\begin{lem}
Let $A,B \in L^\infty$. Then $AH^2_+ \cap B H^2_- = A\, \k^+ S_{A,B}$.
\end{lem}

\beginpf
If $A\phi_+ = B \phi_-$, with $\phi_\pm \in H^2_\pm$, then
$A\phi_+ - B\phi_-=0$ and $\phi_+ \in \k^+ S_{A,B}$, so $AH^2_+ \cap B H^2_- \subseteq \, A\, \k^+ S_{A,B}$.

Conversely, it is clear that $A\, \k^+ S_{A,B} \subseteq AH^2_+$, and from
$A\phi_+ -  B \phi_-=0$ we see that we also have $A\, \k^+ S_{A,B} \subseteq BH^2_-$.
\endpf

\begin{cor}\label{cor:12jul5.8}
If $\k^+ S_{\alpha B_{1-},B_{2+}}=\{0\}$ then $\GS=\{0\}$.
\end{cor}

Now let $\ol{B_{1-}} = I_1 O_1$ and $B_{2+}=I_2 O_2$ be
inner--outer factorizations, with $I_1,I_2$ inner and $O_1,O_2$ outer. Let $I=I_1I_2$.

\begin{cor} 
If $O_2 \in \G H^\infty$ then
\beq\label{eq:12jul5.21}
\GS= \{0\} \quad \hbox{if}\q \ker T_{\alpha \ol I}=\{0\},
\eeq
i.e., if $\alpha $ ``annihilates'' $K_I$ \cite[Def.~6.1]{CMP16}.
\end{cor}

\beginpf
This is an immediate consequence of Corollary~\ref{cor:12jul5.8} and \eqref{eq:12jul4.18}.
\endpf
With the same notation, from the previous result and Proposition \ref{prop:12jul4.14} we thus have:

\begin{prop}
If $I_1$ and $I_2$ are finite Blaschke products and $O_2 \in \G H^\infty$, then $\GS =\{0\}$ if
\beq
\dim K_\alpha \ge \dim K_{I_1} + \dim K_{I_2}.
\eeq
\end{prop}
Now assume that $\GS = \{0\}$. In that case, by Theorem~\ref{thm:12jul5.5} we
have \eqref{eq:12jul5.16}. To characterise $\GK$, we will use the following result, with the notation
\beq
\GH_X = \{ h \in \GH: Xh \in \GH \}.
\eeq

\begin{prop}\label{12jul5.11}
Let $A_{2+},B_{2+} \in  H^\infty$ with
$\inf_{z \in \DD} |A_{2+}(z)|+|B_{2+}(z)|>0$. Then
\beq\label{12jul5.24}
(H^2_+)_{A_{2+}/B_{2+}} = \ker H_{A_{2+}/B_{2+}}= B_{2+}H^2_+.
\eeq
Analogously, if 
$A_{1-},B_{1-} \in  \ol{H^\infty}$ with
$\inf_{z \in \DD} |\ol{A_{1-}}(z)|+|\ol{B_{1-}}(z)|>0$, then
\beq\label{12jul5.25}
(H^2_-)_{A_{1-}/B_{1-}} = \ker \wt H_{A_{1-}/B_{1-}} = B_{1-} H^2_-.
\eeq
\end{prop}

\beginpf
Clearly, $B_{2+}H^2_+ \subseteq (H^2_+)_{A_{2+}/B_{2+}} $. Conversely, if
$\dfrac{A_{2+}}{B_{2+}} \phi_{1+} = \psi_+$, with $\phi_{1+},\psi_+ \in H^2_+$, and $A_{2+}\wt A_{2+}+B_{2+}\wt B_{2+}=1$, then
\begin{flalign*}
\wt A_{2+}  A_{2+} \phi_{1+} = \wt A_{2+} B_{2+} \psi_+ 
&\iff \phi_{1+}-B_{2+}\wt B_{2_+} \phi_{1+} = \wt A_{2+}B_{2+} \psi_+\\
&\iff
\phi_{1+}= B_{2+} (\wt A_{2+}\psi_+ + \wt B_{2+}\phi_{1+}) \in B_{2+}H^2_+.
\end{flalign*}
The proof of \eqref{12jul5.25} is analogous.
\endpf

\begin{rem}{\rm
Note   in \eqref{12jul5.24} that $(H^2_+)_{A_{2+}/B_{2+}}$ can be seen as the kernel of
a possibly unbounded Hankel operator, and analogously for \eqref{12jul5.25}.
Thus Proposition~\ref{12jul5.11} establishes a link between the corona theorem and kernels of Hankel operators
with symbols of a certain form.
}
\end{rem}

Thus we have, using Proposition \ref{12jul5.11}:

\begin{prop}
Let $f,g,\til f,\til g$ be as defined in \eqref{eq:11jul5.12} and \eqref{eq:11jul5.13}. Then
\beq
\ker H_{f \til f^T}=(B_{2+}H^2_+,H^2_+)
\eeq
and
\beq
\ker \wt H_{g\til g^T}=(B_{1-}H^2_-,H^2_-).
\eeq
\end{prop}
Finally, from Theorem~\ref{thm:12jul5.5} we get the following.

\begin{thm}\label{thm:12jul5.14}
Let $\phi$ be defined as in \eqref{eq:11jul5.1} and suppose that $\GS=\{0\}$. Then
$\ker A^{\th,\alpha}_\phi$ does not depend on $\alpha$ and we have
\beq
\ker A^{\th,\alpha}_\phi = B_{2+} \k^+ S_{\ol\th B_{2+},-B_{1-}}.
\eeq
Further, if $\ol B_{1-}= I_1 O_1$ and $B_{2+}=I_2 O_2$ are inner--outer factorizations,
with $O_1,O_2 \in \G H^\infty$, then
\beq\label{eq:12jul5.29}
\ker A^{\th,\alpha}_\phi  = I_2 \ker T_{\ol\th (I_1 I_2)}.
\eeq
\end{thm}

\beginpf
The first part follows from Theorem \ref{thm:12jul5.5} and \eqref{eq:12jul5.29}
follows from Propositions \ref{prop:9jul4.3} and \ref{prop:9jul4.7}.
\endpf

Recall that the question whether Toeplitz kernels of the form arising in \eqref{eq:12jul5.29} are 
equal to $\{0\}$ was discussed in Section~\ref{sec:3}.
We also note that the results here do not depend on $A_{1-}$
nor $A_{2+}$.\\

Finally, let us consider the particular case of ATTO of finite rank, studied in \cite{CP24}, where the symbol takes the form
\beq\label{eq:12julI}
\phi = \ol\th R_+ - \alpha R_- + \sum_{j=1}^N 
\frac{\ol\th P^\alpha_{n_j-1}(z)-\alpha P^{\ol\th}_{n_j-1}(t_j)}
{(z-t_j)^{n_j}},
\eeq
where
\begin{enumerate}[(i)]
\item $R_\pm$ are rational functions vanishing at $\infty$, such that
$R_+ \in H^\infty$ and $R_- \in \ol{H^\infty}$;
\item $t_j \in \TT$ for $j=1,2,\ldots,N$, and are regular points for $\theta$ and $\alpha$;
\item $P^\alpha_{n_j-1}$ and $P^{\ol\th}_{n_j-1}$ are the Taylor polynomials of order $n_j-1$ relative to the point $t_j$, for $\alpha$ and $\ol\th$ respectively.
\end{enumerate}
The formula in \eqref{eq:12julI} can be rewritten in the form
\beq
\phi = \ol\th \frac{Q_1}{\GE D_{1+}}-\alpha \frac{Q_2}{\GE D_{2-}} \in L^\infty,
\eeq
where:
\begin{itemize}
\item $Q_1$ and $Q_2$ are polynomials, sharing no common zeros with $\GE D_{1+}$ and
$\GE D_{2-}$, respectively;
\item $D_{1+}$ is a polynomial of degree $N_1$ with all its zeros in $\DD^e:=\{z \in \CC: |z|>1 \}$;
 \item $D_{2-}$ is a polynomial of degree $N_2$ with all its zeros in $\DD$;
 \item $\GE$ is a polynomial of degree $N_0$ with all its zeros on $\TT$;
 \item $\deg Q_1 < N_0+N_1$ and $\deg Q_2 < N_0+N_2$.
 \end{itemize}
 Therefore
 \beq
 \phi= \ol\th \frac{A_{1-}}{B_{1-}} - \alpha \frac{A_{2+}}{B_{2+}},
 \eeq
 with
 \[
 \begin{array}{ll}
 A_{1-}= \frac{Q_1}{z^{N_0+N_1}} & B_{1-}= \frac{\GE D_{1+}}{z^{N_0+N_1}} \\ \\
 A_{2+}=Q_2 & B_{2_+} = \GE D_{2-}.
 \end{array}
 \]
 By Corollary \ref{cor:12jul5.8} we know that $\GS=\{0\}$ if
 $\ds \k^+ S_{\alpha \GE D_{1+}/z^{N_0+N_1},\GE D_{2-}}=\{0\}$.
 Now,
 \begin{flalign*}
 \k^+ S_{\alpha\GE D_{1+}/z^{N_0+N_1},\GE D_{2-}} &=
 \k^+ S_{\alpha  D_{1+}/z^{N_0+N_1},  D_{2-}}\\
 &= \k^+ S_{\alpha  \ol z^{N_0+N_1}D_{1+},z^{N_2} (D_{2-}/z^{N_2})}\\
 &= D_{1+}^{-1} \k^+ S_{\alpha \ol z^{N_0+N_1},z^{N_2}} \\
 &= D_{1+}^{-1} \ker T_{\alpha \ol z^{N_0+N_1+N_2}}
 \end{flalign*}
 because $D_{1+} \in \G H^\infty$, $D_{2-}/z^{N_2} \in \G \ol{H^\infty}$
 (by Propositions \ref{prop:9jul4.3} and \ref{prop:9jul4.7} and by \eqref{eq:12jul4.6}). Therefore
 \[
 \GS=\{0\} \qquad \hbox{if} \q \deg\alpha > N_0+N_1+N_2.
 \]
Now, by Theorem \ref{thm:12jul5.14},
\begin{flalign*}
\ker A_\phi^{\th,\alpha} &= \GE D_{2-} \k^+ S_{\ol\th \GE D_{2-},-\GE D_{1+}/z^{N_0+N_1}}\\
&= \GE D_{2-} \k^+ S_{\ol\th   D_{2-},-  D_{1+}/z^{N_0+N_1}}\\
&= \GE D_{2-} \k^+ S_{\ol\th z^{N_2} (D_{2-}/z^{N_2})z^{N_0+N_1},-D_{1+}}\\
&= \GE D_{2-} D_{1+} \k^+ S_{\ol\th z^{N_0+N_1+N_2},-1}\\
&= \GE D_{2-} D_{1+} \ker T_{\ol\th z^{N_0+N_1+N_2}}\,,
\end{flalign*}
by Propositions \ref{prop:9jul4.3} and \ref{prop:9jul4.7}, as above. Thus we obtain the same result as in
\cite{CP24}, but in a simpler way.

\subsection*{Acknowledgment}
Work  partially funded by FCT/Portugal through project UIDB/04459/2020 with DOI identifier 10-54499/UIDP/04459/2020

\end{document}